\documentclass{llncs}

\usepackage{hyperref}
\usepackage{graphicx}
\usepackage{stmaryrd}
\usepackage{subfig}
\usepackage{latexsym}
\usepackage{amsmath}
\usepackage{amssymb}
\usepackage{inputenc}
\usepackage{url}
\usepackage{pstricks}
\usepackage{tabularx}

\newcommand{\pc}{\mathbf{P}}
\newcommand{\oc}{\mathbf{O}}
\newcommand{\uc}{\mathbf{U}}
\newcommand{\pp}{\mathbb{P}}
\newcommand{\ox}{\mathbb{X}}
\newcommand{\po}{\mathbb{O}}
\newcommand{\ml}{\mathfrak{L}}
\newcommand{\mm}{\blacksquare}

\begin{document}

\title{A Program in Dialectical Rough Set Theory}
\author{A. Mani}
\institute{ Calcutta Mathematical Society\\
9/1B, Jatin Bagchi Road\\
Kolkata(Calcutta)-700029, India\\
E-Mail: \texttt{$a.mani@member.ams.org$}\\
Homepage: \url{http://www.logicamani.co.cc}}

\maketitle

\begin{abstract}
A dialectical rough set theory focussed on the relation between roughly equivalent objects and classical objects was introduced in \cite{AM699} by the present author. The focus of our investigation is on elucidating the minimal conditions on the nature of granularity, underlying semantic domain and nature of the general rough set theories (RST) involved for possible extension of the semantics to more general RST on a paradigm. On this basis we also formulate a program in dialectical rough set theory. The dialectical approach provides better semantics in many difficult cases and helps in formalizing a wide variety of concepts and notions that remain untamed at meta levels in the usual approaches. This is a brief version of a more detailed forthcoming paper by the present author. 
\end{abstract}

\section{Introduction}
In the present author's formalism a binary predicate or a predicate involving two types of objects that are in some sense opposed to each other in their properties is a \emph{dialectical predicate}. A theory focussing on one or more dialectical relationships and making use of associated dialectical rules is a \emph{dialectical theory}. The concept of dialectical predicate may also be required to be quasi-inductive by way of requiring the existence of an associated dialectical theory and so they may have complex external meaning. By a \emph{dialectical RST} we understand a rough set theory imbibed or extended with dialectical semantic features.

It is well known that different types of semantics for RST can be developed with different meta level assumptions. The latter often lead to some loss or distortion of information, externalization of concepts and measures, and provides for what may be called 'point of view semantics'. For example the Lesniewski ontology based mereological approaches \cite{PLS} and the approach in this paper do not share many meta-level concepts. The process of combining such semantics (which can be seen as the next level of progression) leads to dialectical relationships between different types of objects and processes. In \cite{AM699}, in the context of classical RST, we focus on the dialectical relationship between \emph{roughly equivalent objects} in the rough semantic domain and objects in the classical semantic domain. Of the two dialectical semantics developed, one uses a latent dialectical predicate, while the other makes it explicit. Various other types of algebraic or partial algebraic rough semantics (like \cite{AM105}, \cite{AM3}) can be \emph{fibred} under restrictions on dialectical relationships along similar lines. But our focus will be on the essential parts of the afore-mentioned dialectical semantics for prospective generalization.

The focus of the present paper is relevant for many contexts including those of dynamic spaces or multiple approximation spaces \cite{PP4}, \cite{DO} and mathematical morphology. It is in itself a general approach to multi source approximation spaces and other less tractable general Approximation Spaces ($\mathcal{AS}$). The exact realisation of general rough set theories as special semi-set theories seems possible through dialectical RST. The general motivations for dialectical RST are also to provide a framework for investigating relative distortions introduced by different theories, to improve the interface between rough and classical semantic domains in application contexts, to address issues relating to truth at the semantic level and to develop dialectical logics for rough semantics.  

\section{Semantic Domains, Object and Meta Levels}

The term \emph{object level} will mean a description that can be used to directly interface with fragments (sufficient for the theories or observations under consideration) of the concrete real world. Meta levels concern fragments of theories that address aspects of dynamics at lower meta levels or the object level. Importantly, we permit meta level aspects to filter down to object levels relative a different object level specification. So we can always construct more meta levels and expressions constructed by us in these carry intentions (in the modern sense). 

A semantic domain $\mathcal{C}$ (or a domain of discourse) can be seen as a function from some types to a category of generalized or multi-sorted relational structures ( generalized structures of different types are taken as objects and the morphisms are morphisms of generalized structures) that has the following properties: functor, closure, richness, finite reduct, renaming, pair and diagram (see \cite{MD}). When many of the properties of possible semantic models in ZFC, for example are known, then it is sensible to speak of semantic domains determined by objects sharing certain properties. Of course, there are no real methodological hurdles in this.

In classical rough set theory, we start with an $\mathcal{AS}$ $\left\langle S,\,R \right\rangle $, with $R$ being an equivalence on the set $S$. On the power set $\wp (S)$, we can define lower and upper approximation operators apart from the usual Boolean operations. This constitutes a semantics for rough set theory (though not satisfactory) in a classicalist perspective. This continues to be true even when $R$ is some other type of binary relation. For more general purposes, we can replace $\wp (S)$ by some set with a parthood relation and some approximation operators defined on it. The associated semantic domain in the sense of a collection of restrictions on possible objects, predicates, constants, functions and low level operations on those will be referred to as the classical semantic domain for general RST. In contrast the semantic domain associated with sets of roughly equivalent or relatively indiscernible objects with respect to this domain will be called the rough semantic domain. Though we can generate many other semantic domains including hybrid ones for different types of rough semantics, these two broad domains will always be. 

Formal versions of these types of semantic domains will be useful for clarifying the relation with categorical approaches to fragments of granular computing as considered in \cite{BY}. But even without a formal definition it can be seen that two approaches are not equivalent. A major difference is the requirement of an a-priori description of granularity in the categorical approach. This is not a desirable state of affairs as the most suitable form of granules should evolve relative particular semantics or semantic domains. The entire category \textbf{ROUGH} of rough sets in \cite{BY}, assumes a uniform semantic domain as evidenced by the notion of objects and morphisms. A unifying semantic domain may not be definable for many sets of semantic domains in our approach. This tells us that the categorical approach needs to be extended to provide an accessible setting. 

\section{General Rough $Y$-Systems}
Different formal definitions of granules have been used in the literature on rough sets and in granular computing. We propose a refined version capable of handling most contexts and permitting relaxation of set-theoretic axioms at a later stage. Towards this we introduce a generalized form of \emph{rough orders} (\cite{IT2}) to include relation based RST, cover based RST and more. These structures are provided with enough structure so that we can associate a classical semantic domain and at least one rough semantic domain of roughly equivalent objects along with admissible operations and predicates. 
\begin{definition}
A \emph{general rough Y system} (\textsf{RYS}) will be a tuple of the form  $\left\langle {S},\,W,\,\pc ,\,(l_{i})_{1}^{n},\,(u_{i})_{1}^{n},\,+,\,\cdot,\,\sim,\,1 \right\rangle $ satisfying all of the following ($\pc$ is intended as a binary relation on $S$ and $W\,\subset\,S$, $n$ being a finite positive integer. $\iota$ is the description operator of FOPL: $\iota(x) \Phi(x)$ means 'the $x$ such that $\Phi(x)$ ' ):
\begin{enumerate}
\item {$(\forall x) \pc xx$ ; $(\forall x,\,y)(\pc xy,\,\pc yx\,\longrightarrow\,x=y)$}
\item {For each $i,\,j$, $l_{i}$, $u_{j}$ are surjective functions $:S\,\longmapsto\,W$ }
\item {For each $i$, $(\forall{x,\,y})(\pc xy\,\longrightarrow\,\pc (l_{i}x)(l_{i}y))$}
\item {For each $j$, $(\forall{x,\,y})(\pc xy\,\longrightarrow\,\pc (u_{i}x)(u_{i}y))$}
\item {For each $i$, $(\forall{x})\,\pc (l_{i}x)x,\,\pc (x)(u_{i}x)$}
\item {For each $i$, $(\forall{x})(\pc (u_{i}x)(l_{i}x)\,\longrightarrow\,x=l_{i}x=u_{i}x)$}
\item {For each $i$, $(\forall{x})\,\pc (l_{i}x)(u_{i}l_{i}x),\,\pc (l_{i}u_{i}x)(u_{i}x)$}
\end{enumerate}

The operations $+,\,\cdot$ and the derived operations $\oc,\, \pp,\,\uc,\,\ox,\,\po $ will be defined as follows:
\begin{description}
\item [Overlap: ]{$\oc xy\,\mathrm{iff} \,(\exists z)\,\pc zx\,\wedge\,\pc zy$}
\item [Underlap: ]{$\uc xy\,\mathrm{iff} \,(\exists z)\,\pc xz\,\wedge\,\pc yz$} 
\item [Proper Part: ]{$\pp xy\,\mathrm{iff} \,\pc xy\,\wedge\,\neg \pc yx$}
\item [Overcrossing: ]{$\ox xy \,\mathrm{iff}\,\oc xy\,\wedge\,\neg \pc xy$}
\item [Proper Overlap: ]{$\po xy \,\mathrm{iff}\,\ox xy \,\wedge\,\ox yx $}
\item[Sum: ]{$x+y=\iota z (\forall w)(\oc wz\, \leftrightarrow\,(\oc wx \vee \oc wy))$}
\item[Product: ]{$x \cdot y=\iota z (\forall w)(\pc wz\,\leftrightarrow\,(\pc wx \wedge \pc wy))$ }
\item[Difference: ]{$x - y=\iota z (\forall w)(\pc wz\,\leftrightarrow\,(\pc wx \wedge \neg \oc wy))$}
\item[Associativity:]{We will assume that $+,\,\cdot$ are associative operations.}
\end{description}
\end{definition}

A \textsf{RYS} is intended to capture a minimal common fragment of different rough set theories. Note that the parthood relation $\pc$, taken as a general form of rough inclusion (in a suitable semantic domain), is not necessarily transitive. Transitivity of $\pc$ is a sign of fair choice of attributes (at that level of classification), but non transitivity may result by the addition of attributes and so the property by itself says little about the quality of classification. Examples (better than the \emph{handle-door-house} example) are easy to construct. In classical rough set theory, 'supplementation' in the stricter sense (that is, $(\neg \pc xy\,\longrightarrow\,\exists z (\pc zx\,\wedge\,\neg \po zy))$ does not hold), while the weaker version $(\neg \pc xy\,\longrightarrow\,\exists z (\pc zx\,\wedge\,\neg \oc zy))$ is trivially satisfied due to the existence of the empty object ($\emptyset$).  Proper selection of semantic domains is essential for avoiding possible problems of ontological innocence (\cite{MCE}), wherein the 'sum' operation may result in non existent objects relative the domain. A similar operation results in 'plural reference' in \cite{PLS} and related papers. 

\begin{definition}
In the above, two approximation operators $u_i$ and $l_i$ will be said to be {\em dual \/} to each other if and only if $(\forall A\subset S)\,A^{u_{i}l_{i}}=A^{u_{i}} \;;\; A^{l_{i} u_{i}}=A^{l_{i}}$. 
\end{definition}

Possible defining properties of a set of granules include all of the following: ($t_{i},\,s_{i}$ are term functions formed with $+,\,\cdot,\,\sim$, while $p,\,r$ are finite positive integers. $\mathbf{\forall}i$, $\mathbf{\exists}i$ are meta level abbreviations)

\begin{description}
\item [Representation of Approximations, RA] {$\mathbf{\forall}i$, $(\forall x)(\exists y_{1},\ldots y_{r}\in \mathcal{G})$ $y_{1} + y_{2} + \ldots + y_{r}=x^{l_{i}}$ and $(\forall x)(\exists y_{1},\,\ldots\,y_{p}\in \mathcal{G})\,y_{1} + y_{2} + \ldots + y_{p}=x^{u_{i}}$}
\item [Weak Representation of Approximations, WRA] {$\mathbf{\forall}i$, $(\forall x \exists y_{1},\ldots y_{r}\in \mathcal{G})$ $t_{i}(y_{1},\,y_{2}, \ldots \,y_{r})=x^{l_{i}}$ and $(\forall x)(\exists y_{1},\,\ldots\,y_{r}\in \mathcal{G})\,t_{i}(y_{1},\,y_{2}, \ldots \,y_{p}) = x^{u_{i}}$}
\item [Absolute Crispness of Granules, ACG] {For each $i$, $(\forall y\in \mathcal{G})\,y^{l_{i}} = y^{u_{i}} = y$}
\item [Weak Crispness of Granules,WCG] {$\mathbf{\exists}i$, $(\forall y\in \mathcal{G}) y^{l_{i}} = y^{u_{i}} = y$.}
\item [Mereological Atomicity,MER] {$\mathbf{\exists}i$, $(\forall y\in \mathcal{G})(\forall x \in S)(\pc xy,\,x^{l_{i}} = x^{u_{i}} = x \longrightarrow x = y)$}
\item [Lower Stability,LS] {$\mathbf{\exists}i$, $(\forall y \in \mathcal{G})(\forall {x\in S})\,(\pc yx\,\longrightarrow\,\pc (y)(x^{l_{i}}))$}
\item [Upper Stability,US] {$\mathbf{\exists}i$, $(\forall y\in\mathcal{G})(\forall {x\in S})\,(\oc yx \longrightarrow P(y)(x^{u_{i}}))$}
\item [Stability, ST] {Shall be the same as the satisfaction of lower and upper stability}
\item [Absolute Stability,AS] {Is the same as the satisfaction of stability for every $i$}
\item [No Overlap, NO] {$(\forall x ,\,y\in\mathcal{G}) \neg \po xy$}
\item [Full Underlap, FU] { $\mathbf{\exists}i$, $(\forall x,\,y\in\mathcal{G})(\exists z\in S ) \pp xz,\,\pp yz,\,z^{l_{i}} = z^{u_{i}} = z$}
\item [Unique Underlap, UU] {For at least one $i$, $(\forall x,\,y\in\mathcal{G})(\pp xz,\,\pp yz,\, z^{l_{i}}=z^{u_{i}}=z,\,\pp xb,\,\pp yb,\, b^{l_{i}}=b^{u_{i}}=b\,\longrightarrow\,z=b)$}
\end{description}

\begin{definition}
A subset $\mathcal{G}$ of $S$ in a \textsf{RYS} will be said to be an \emph{admissible set of granules} provided the properties \textsf{WRA}, \textsf{LS} and \textsf{FU} are satisfied by it. 
\end{definition}

In cover based rough set theories (see \cite{AM960} for example) we define different approximations with the help of a determinate collection of subsets. These subsets satisfy the properties \textsf{WRA}, \textsf{LS} and \textsf{FU} and are therefore admissible granules.  But they do not in general have many of the nicer properties of granules in relation to the approximations. However at a later stage we may be able to refine these and construct a better set of granules (see \cite{AM960} for details) for the same approximations. Similar process of refinement can be used in other types of rough set theories as well. For these reasons we shall refer to the former as \emph{initial granules} and the latter as relatively \emph{refined granules}. In the extended example in the next section, these may be seen to be in the 'simplest non-trivial state'.

Note that we are speaking of better set of granules for the same approximations and not about modifications of granules that may result in better approximations on the same strategy of forming approximations. But all of the cases will be relevant in this paper.

\section{Dialectical REQ Rough Set Theory}

By \emph{Dialectical REQ Rough Set Theory} (\textsf{QRST}) we mean an integrated theory based over \textsf{RYS} including multiple types of rough approximations, multiple types of roughly equivalent objects that seeks to study the dialectical relationship between these types of roughly equivalent objects with the help of additional new hybrid operations. Ideally we should expect all of the component theories to be deducible from the \textsf{QRST}. The hybrid operations mentioned may turn out to be extensions of operations on the component theories. Specification of the different types of granules required is more important for these hybrid operations and for easing the handling and interpretation of dialectical predicates at both the semantic level and at the proof-theoretic level. Our reference system for the formulation of the program will be the specific \textsf{QRST}s (\cite{AM699}) below: 

\subsection{Extended Example: QRST-1,2}

Let $S=\left\langle \underline{S},\,R\right\rangle $ be an $\mathcal{AS}$. If $A\,\subset S$, $A^{l}=\bigcup\{[x]\,;\,[x]\,\subseteq\,A\}$ and $A^{u}=\bigcup\{[x]\,;\,[x]\,\cap\,A\,\neq\,\emptyset\}$ are the \emph{lower} and \emph{upper approximation} of $A$ respectively. If $A,\,B\in\wp (S)$, then they are said to be \emph{roughly equal} if and only if $A^{l}=B^{l}$ and $A^{u}=B^{u}$. In this case we write $A\,\approx\,B$. $[x]$ is the equivalence class (w.r.t $\approx$) formed by a $x\in\wp (S)$. Our base set will be $\wp (S)\,\cup\,(\wp (S)|\approx)$ as opposed to $\wp (S)|\approx$ (used in the construction of a pre-rough set algebra, \cite{BC1}). The new operations $\oplus, \, \odot$ introduced below permits us to combine objects (form general unions and intersections respectively) in the rough and classical domain to form objects in the rough domain. This is not possible in classical RST, where classicalist object level elements cannot interact with objects in the rough semantic domain. Possibly inconsistent interpretations can also be accommodated through the interaction. 

\begin{definition}
On $Y=\wp (S)\,\cup\,\wp (S)|\approx$, the operations $\ml,\,\oplus,\,\odot,\,\mm,\,\rightsquigarrow,\,\twoheadrightarrow,\,\sim$ will be defined as follows: (we assume that the operations $\cup,\,\cap,\,^{c},\,^{l},\,^{u}$ and $\sqcup,\,\sqcap,\,L,\,M,\,\neg,\,\Rightarrow $ are available on $\wp (S)$ and $\wp(S)|\approx$ respectively. Further $\tau_{1} x \Leftrightarrow x\in \wp(S)$ and $\tau_{2}x \Leftrightarrow x\in \wp(S)|\approx$)
\begin{itemize}
\item {$
\ml x = \left\lbrace  \begin{array}{ll}
 x^{l} & \textrm{if $\tau_{1}x$}\\
 Lx & \textrm{if $\tau_{2}x$}\\
 \end{array} \right.$
}
\item {$
\mm x = \left\{ \begin{array}{ll}
 x^{u} & \textrm{if $\tau_{1}x$}\\
 \neg L \neg x & \textrm{if $\tau_{2}x$}\\
 \end{array} \right.
$}
\item {$
 x\,\oplus\,y = \left\{ \begin{array}{ll}
 x\,\cup\,y & \textrm{if $\tau_{1}x,\,\tau_{1}y$}\\
\left[ x\,\cup\,(\bigcup_{z\in y} z)\right]  & \textrm{if $\tau_{1}x,\,\tau_{2}y$}\\
\left[ (\bigcup_{z\in x} z)\,\cup\,y\right]  & \textrm{if $\tau_{2}x,\,\tau_{1}y$} \\
 x\,\sqcup\,y & \textrm{if $\tau_{2}x,\,\tau_{2}y$}
 \end{array} \right.
$} 
\item {$
 x\,\odot\,y = \left\{ \begin{array}{ll}
 x\,\cap\,y & \textrm{if $\tau_{1}x,\,\tau_{1}y$}\\
 \left[ x\,\cap\,(\bigcap_{z\in y} z)\right]  & \textrm{if $\tau_{1}x,\,\tau_{2}y$}\\
 \left[ (\bigcap_{z\in x}\,z)\,\cap\,y\right]  & \textrm{if $\tau_{2}x,\,\tau_{1}y$ }\\
 x\,\sqcap\,y & \textrm{if $\tau_{2}x,\,\tau_{2}y$}
 \end{array} \right.
$}
\item {$
 \sim x\, = \left\{ \begin{array}{ll}
 x^{c} & \textrm{if $\tau_{1}x$}\\
 \neg x & \textrm{if $\tau_{2}x$}
 \end{array} \right.
$}
\item {$
 x\,\rightsquigarrow\,y = \left\{ \begin{array}{ll}
 x\,\cup\,y^{c} & \textrm{if $\tau_{1}x,\,\tau_{1}y$}\\
 \left[\bigcup_{z\in y} (x\,\cup\, z^{c}\right] & \textrm{if $\tau_{1}x,\,\tau_{2}y$}\\
 \left[\bigcup_{z\in x}(z\,\cup\,y^{c}\right] & \textrm{if $\tau_{2}x,\,\tau_{1}y$ }\\
 (\neg{L}(x)\sqcup{L}(y))\sqcap(L(\neg{x})\sqcup\neg{L}(\neg{y})) & \textrm{if $\tau_{2}x,\,\tau_{2}y$}
 \end{array} \right.
$}
\item {$
 x\,\twoheadrightarrow\,y = \left\{ \begin{array}{ll}
 [x\,\cup\,y^{c}] & \textrm{if $\tau_{1}x,\,\tau_{1}y$}\\
 \left[\bigcup_{z\in y} (x\,\cup\, z^{c}\right] & \textrm{if $\tau_{1}x,\,\tau_{2}y$}\\
 \left[\bigcup_{z\in x}(z\,\cup\,y^{c}\right] & \textrm{if $\tau_{2}x,\,\tau_{1}y$ }\\
 (\neg{L}(x)\sqcup{L}(y))\,\sqcap\,(L(\neg{x})\sqcup\neg{L}(\neg{y})) & \textrm{if $\tau_{2}x,\,\tau_{2}y$}
 \end{array} \right.
$} 
\end{itemize}
\end{definition}
\begin{definition}
In the above context a partial algebra of the form\\ $W=\left\langle \underline{\wp (S)\cup \wp(S)|\approx},\,\neg,\,\sim,\,\oplus,\,\odot,\,\mm,\,\ml,\,0,\,1,\,\bot,\,\top (1,\,1,\,2,\,2,\,1,\,1,\,0,\,0,\,0,\,0)\right\rangle$ will be called a \emph{concrete enriched pre-rough algebra} (CERA) if a pre-rough algebra structure is induced on $\wp(S)|\approx$. \emph{Concrete enriched rough algebras} can be defined in the same manner. If the $\mathcal{AS}$ is $X$, then the derived CERA will be denoted by $\mathfrak{W}(X)$. Note that the two implication-like operations are definable in terms of the others. CERAs are well defined because of the representation theory of pre-rough algebras.
\end{definition}
\begin{theorem}
A CERA satisfies all of the following: (The first two conditions state that the $\tau_{i}$s are abbreviations)
\begin{description}
\item [type-1]{$(x\,\rightsquigarrow\,x=\top \,\longleftrightarrow\,\tau_{1} x)$}
\item [type-2]{$(\neg x=\neg x \,\longleftrightarrow\,\tau_{2} x )$}
\item [ov-1]{$\sim \sim x=\sim x$ ; $\ml \ml x=\ml x $ ; $\mm \ml x=\ml x$ }
\item [ov-2]{$\ml x\,\oplus\,x=x $ ; $\ml x\,\odot\,x=\ml x $ ; $\mm x\,\oplus\,x=\mm x $ ; $\mm x\,\odot\,x=x$}
\item [ov-3]{$\ml \mm x=\mm x $ ; $x\,\,\oplus\,x=x $ ; $x\,\odot\,x=x $}
\item [qov-1]{$(\tau_{1}x\,\longrightarrow\,\sim x\,\oplus\,  x=\top)$ ; $(\tau_{2}x\,\longrightarrow\,\sim \ml x\,\oplus\,\ml x=1)$}
\item [qov-2]{$\sim \bot=\top$ ; $\sim 0=1$}
\item [u1]{$x\,\oplus\,(x\,\oplus\,(x\,\oplus\,y))=x\,\oplus\,(x\,\oplus\,y) $ ; $x\,\odot\,(x\,\odot\,(x\,\odot\,y))=x\,\odot\,(x\,\odot\,y) $}
\item [u2]{$x\,\oplus\,y=y\,\oplus\,x$; $x\,\odot\,y=y\,\odot\,x$}
\item [ter(i)]{$(\tau_{i}x, \tau_{i}y, \tau_{i}z\,\longrightarrow\,x\,\oplus\,(y\,\oplus\,z)=(x\,\oplus\,y)\,\oplus\,z,\,\,x\,\oplus\,(y\,\odot\,z)=(x\,\oplus\,y)\odot (x\,\oplus\,z),\,\,x\,\odot\,(y\,\odot\,z)=(x\,\odot\,y)\,\odot\,z)$ ; i=1,2}
\item [bi(i)]{$(\tau_{i}x, \tau_{i}y\,\longrightarrow\,x\,\oplus\,(x\,\odot\,y)=x,\,\,\sim (x\,\odot\,y)=\sim x\,\oplus\,\sim y)$ ; i=1,2}
\item [bm]{$(\tau_{1}x,\,\tau_{2}y,\,x\,\oplus\,y=y\,\longrightarrow\,\mm x\,\oplus\,y=y)$}
\item [hra1]{$(\tau_{1}x,\,(1\,\odot\,x=y)\,\vee\,(y=x\,\oplus\,0)\longrightarrow\,\tau_{2}y) $}
\end{description}
\end{theorem}

\begin{definition}
An \emph{abstract enriched pre-rough partial algebra} (AERA) will be a partial algebra of the form $S=\left\langle\underline{S},\,\neg,\,\sim,\,\oplus,\,\odot,\,\mm,\,\ml,\,0,\,1,\,\bot,\,\top  \right\rangle$ satisfying :
\begin{description}
\item [RA]{$\mathrm{dom}(\neg)$ along with the operations $(\oplus,\,\odot,\,\mm,\,\ml,\,\sim,\,0,\,1)$ restricted to it and the definable $\Rightarrow$ forms a pre-rough algebra}
\item [BA]{$\underline{S}\setminus \mathrm{dom}(\neg)$ with the operations $(\oplus,\,\odot,\,\mm,\,\ml,\,\sim,\,\top,\,\bot)$ restricted to it forms a topological boolean algebra (with an interior and closure operator)}
\item [IN]{Given the definitions type-1, type-2, all of u1, u2, ter(i), bi(i), bm and hra hold}
\end{description}
\end{definition}

\begin{theorem}[Representation \cite{AM699}]
Every AERA $S$ has an associated $\mathcal{AS}$ $X$ (up to isomorphism), such that the derived CERA $\mathfrak{W}(X)$ is isomorphic to it. AERAs are also definable by a set of quasi equations.
\end{theorem}

In the above semantics it is not possible to transform objects in the rough domain to objects in the classical domain. But if we take the universe to be the set of tuples of the form $\{(x,\,0 \oplus x)\,:\,\tau_{1}x\}\,\cup\,\{(y,\,x)\,:\,\tau_{2} y,\,\tau_{1} x,\,x\,\oplus\,0=y\}=K$ ($x,\,y$ being elements of a CERA), the requirement is trivialised. The algebra formed by defining similar extended operations is called a \emph{concrete rough dialectical algebra} (CRAD) and associated dialectical logics can be defined. 

\section{The Dialectical REQ Rough Set Program}
Many partial/total algebraic semantics for \textsf{QRST} can be constructed in the rough semantic domain determined by the choice of definition of rough equalities along the following lines:
\begin{enumerate}
\item {Define a rough equality $\approx$ via $x\,\approx\,y$ iff a boolean combination of equalities of the form $x^p=y^p$ hold for approximation operators $p$.}
\item {Form the quotient $S|\approx$. Clearly there is an injection from $S|\approx$ into $W$ (usually it will be a bijection) and we can define approximation operators $P$ on $S| \approx$ via $P[x]=[x^{p}]$, '$[]$' being the operation of forming quotients.}
\item {Other induced operations corresponding to $+,\,\cdot$ and others may be partial. The partial algebra may admit of good representation theorems. In this case we have a good semantics.}
\item {Otherwise we can use a wide variety of algebraic approaches including the approach used for bitten rough set theory in \cite{AM105} or the choice inclusive approach in \cite{AM69} or others.}
\item {Even if a good algebraic semantics is not possible, the dialectical integration can result in semantics with better representation theory.}
\end{enumerate}

Based on the above possibilities, our proposed \emph{program} consists in forming:
\begin{description}
\item [A]{An integrated algebraic semantics with hybrid operations as in the extended example above or over a universe of ordered pairs of the form $(a,b)$ with each being of different type and being dialectically opposed to the other. This is also relevant when we try to use different types of granules for constructing structurally similar approximations (and their consequent dialectics is of interest, \cite{AM960}).}
\item [B]{Given the semantics, our next task would be prove both abstract and concrete representation theorems and possibly internalise membership functions.}
\item [C]{Sequent calculi taking the dialectical aspect into consideration would be our next aim. In most cases it is possible to use a set of dialectical, non-dialectical rules and interface rules along with typed variables and constants (see \cite{AM933}).}
\end{description}

Relaxation of the ZF axioms used in the definition of \textsf{RYS} is bound to have difficult consequences as the relevant duality based representation theory is not well developed.

\bibliographystyle{splncs.bst}
\bibliography{newsem36.bib}

\begin{thebibliography}{10}

\bibitem{AM699}
Mani, A.:
\newblock Integrated dialectical logics for relativised general rough set
  theory.
\newblock In: Internat. Conf. on Rough Sets, Fuzzy Sets and Soft Computing,
  Agartala, India. (2009)  6pp

\bibitem{PLS}
Polkowski, L., Polkowska, S.M.:
\newblock Reasoning about concepts by rough mereological logics.
\newblock In Wang, G.,  et~al., eds.: RSKT 2008, LNAI 5009, Springer (2008)
  197--204

\bibitem{AM105}
Mani, A.:
\newblock Algebraic semantics of similarity-based bitten rough set theory.
\newblock Fundamenta Informaticae (2009)  1000--1021

\bibitem{AM3}
Mani, A.:
\newblock Super rough semantics.
\newblock Fundamenta Informaticae \textbf{65} (2005)  249--261

\bibitem{PP4}
Pagliani, P.:
\newblock Pretopologies and dynamic spaces.
\newblock Fundamenta Informaticae \textbf{59} (2004)  221--239

\bibitem{DO}
Demri, S., Orlowska, E.:
\newblock Incomplete Information: Structures, Inference, Complexity.
\newblock Springer-Verlag (2002)

\bibitem{MD}
Mundici, D.:
\newblock Generalization of abstract model theory.
\newblock Fundamenta Math \textbf{124} (1984)  1--25

\bibitem{BY}
Banerjee, M., Yao, Y.:
\newblock Categorical basis for granular computing.
\newblock In Stefanowski, J., Ramanna, S., Butz, C.J., Pedrycz, W., Wang, G.,
  eds.: RSFDGrC 2007, LNAI.
\newblock Springer (2007)  427--434

\bibitem{IT2}
Iwinski, T.B.:
\newblock Rough orders and rough concepts.
\newblock Bull. Pol. Acad. Sci (Math) \textbf{(3--4)} (1988)  187--192

\bibitem{MCE}
Carrara, M., Martino, E.:
\newblock On the ontological commitment of mereology.
\newblock Rev. Symb. Logic \textbf{2} (2009)  164--174

\bibitem{AM960}
Mani, A.:
\newblock Towards an algebraic approach for cover based rough semantics and
  combinations of approximation spaces.
\newblock In: RSFDGrC 2009, LNAI. (2009)  8pp

\bibitem{BC1}
Banerjee, M., Chakraborty, M.K.:
\newblock Rough sets through algebraic logic.
\newblock Fundamenta Informaticae \textbf{28} (1996)  211--221

\bibitem{AM69}
Mani, A.:
\newblock Meaning, choice and similarity based rough set theory.
\newblock Internat. Conf. Logic and Applications, Chennai;
  http://arxiv.org/abs/0905.1352 (2009)  1--12

\bibitem{AM933}
Mani, A.:
\newblock Dialectically presentable logics with applications to fom.
\newblock Preprint (2008)

\end{thebibliography}

\end{document}